\documentstyle[11pt]{article}
\begin{document}

\begin{center}
{\Large \bf Characteristic foliation on the discriminantal
hypersurface of a holomorphic Lagrangian fibration}

\bigskip
 {\large \bf Jun-Muk Hwang}\footnote{This work was supported by the Korea Research Foundation Grant
funded by the Korean Government(MOEHRD)(KRF-2006-341-C00004).}
{\large \bf and Keiji Oguiso}\footnote{This work was supported by JSPS.}

 \end{center}

\bigskip
\begin{abstract}
We give a Kodaira-type classification of general singular fibers of
a holomorphic Lagrangian fibration in Fujiki's class $\mathcal C$.
Our approach is based on the study of the characteristic vector
field of the discriminantal hypersurface, which naturally arises
from the defining equation of the hypersurface via the symplectic
form. As an application, we show that the characteristic foliation
of the discriminantal hypersurface has algebraic leaves which are
either rational curves or smooth elliptic curves.
\end{abstract}

\bigskip
 Key words: holomorphic Lagrangian fibration, characteristic foliation

\medskip
2000MSC:14D06, 14J40, 32M05, 32M25.

\bigskip
\section{Introduction}

We work in the category of complex analytic sets. Let $(M,
\omega)$ be a holomorphic symplectic manifold and $Y \subset M$ be
a reduced hypersurface. We will not assume that $M$ is compact.
The restriction of $\omega$ on the smooth locus $Y_{\rm reg}$ of
$Y$ has a kernel of rank 1 at every point, defining a foliation of
rank 1 on $Y_{\rm reg}$. We will call it the characteristic
foliation on $Y$. (See Section 2 for details.) A leaf of this
foliation is said to be {\it algebraic} if its closure in $M$ is
an algebraic curve, i.e., a compact complex variety of dimension
1. A natural question is the following.

\medskip
{\bf Question 1.1} {\it Are the leaves of the characteristic
foliation algebraic?}

\medskip
The answer is no in general. Note
that if $Y$ is a smooth projective variety with Picard number 1,
it does not admit a foliation with algebraic leaves. For example,
if $M$ is a projective holomorphic symplectic manifold with Picard
number 1 of dimension $\geq 4$, which certainly exists (see e.g.,
[Og]), then for
any smooth hypersurface
$Y$ in $M$, the general leaves of the characteristic foliation are
not algebraic. The essential point of Question 1.1 is to {\it find}
a suitable class of hypersurfaces under which the leaves of the
characteristic foliations are closed.

\medskip
The motivation for the current paper comes from [HR], where the
Hitchin system $f: M \longrightarrow B$ was studied. Here, $M$ is
the moduli of Higgs bundles over a curve and $B$ is a certain
affine space. In [HR], the discriminantal hypersurface $Y \subset
M$, i.e., the locus of singular fibers of $f$,  was studied and it
was noticed ([HR,Remark 4.9]) that the closure of the leaves of
the characteristic foliation on $Y$ are rational curves called
Hecke curves.  In other words, the discriminantal hypersurface of
the Hitchin system was one example where we have a positive answer
for Question 1.1. One of our main results is a generalization of
this fact in the following way. See Sections 2 and 3 for the
precise definitions of a holomorphic Lagrangian fibration and a
discriminantal hypersurface.

\medskip
{\bf Theorem 1.2} {\it Let $M$ be a holomorphic symplectic
manifold and $f: M \longrightarrow B$ be a holomorphic Lagrangian
fibration over a complex manifold $B$. Assume that each fiber of
$f$ is of class ${\cal C}$, i.e., bimeromorphic to a compact
K\"ahler manifold.  Let $Y$ be the discriminantal hypersurface of
$f$. Then the characteristic foliation on $Y$ has algebraic leaves
and  the closures of the leaves are either rational curves or
elliptic curves.  }

\medskip
In the course of proving Theorem 1.2, we are naturally led to
study the structure of general singular fibers of $f$, which will
be of its own interest. The results can be summarized as the
following two theorems (Theorems 1.3 and 1.4). Both theorems
are directly motivated by the
study of the characteristic foliation.

\medskip
{\bf Theorem 1.3} {\it  Let $M$ be a holomorphic symplectic
manifold of dimension $2n$ and $f: M \longrightarrow B$ be a
holomorphic Lagrangian fibration over a complex manifold $B$ of
dimension $n$.  Assume that each fiber of $f$ is of class ${\cal
C}$. Let $\Delta \subset B$ be the hypersurface consisting of
critical values of $f$. For a general point $b$ of an irreducible
component $D$ of $\Delta$, let $X$ be the underlying variety of an
irreducible component of the singular fiber $M_b = f^{-1}(b)$ and
$\hat{X}$ be the normalization of $X$. Then $\hat{X}$ is a compact
complex manifold of class ${\cal C}$ such that the Albanese
variety ${\rm Alb}\,(\hat{X})$ has dimension $n-1$ and the Albanese
map $\alpha: \hat{X} \longrightarrow {\rm Alb}\,(\hat{X})$ is either a
${\mathbf P}_1$-bundle or an elliptic fiber bundle. The latter case
occurs only when $M_b$ is irreducible and non-reduced, in which
case, the reduction $Y_{b} := (M_b)_{\rm red}$ of $M_b$ is smooth. }

\medskip
Theorem 1.3 describes the structure of each irreducible component of
a general singular fiber. The next theorem describes the structure
of the whole general singular fiber.

For the statement, we need a few more terminologies. Under the same notations as in Theorem
1.3, let $Y_b$ be the reduction of each general singular fiber $M_b$. We call
an irreducible curve $\Theta$ on $Y_b$ a {\it characteristic curve} if
it is the image of some Albanese fiber $\alpha : \hat{X} \longrightarrow {\rm Alb}\,(\hat{X})$ of some irreducible component $X$ of $Y_b$. We define two points $y_1$ and $y_2$ on $Y_b$ are {\it equivalent}
if there exist finitely many characteristic curves $\Theta_1, \ldots,
\Theta_N,$ such that $y_1 \in \Theta_1, y_2 \in \Theta_N$ and $\Theta_i
\cap \Theta_{i+1} \neq \emptyset$ for each $1 \leq i \leq N-1$. Then, each
equivalence class is of the form $\cup_{s \in \Lambda} \Theta_s$,
where $\Theta_s$ are characteristic curves. Here the index set $\Lambda$ is
possibly an infinite set.
For each characteristic curve $\Theta_s$, we define the multiplicity $r_s$
to be the multiplicity of the unique irreducible
component $H$ of the discriminant hypersurface $f^{-1}(D)$ such that
$\Theta_s \subset H$. We call the $1$-cycles
$$\sum_{s \in \Lambda} r_s\Theta_s\,\, ,\,\, \sum_{s \in \Lambda}
\Theta_s\,\, ,$$
possibly with infinitely many irreducible components, a {\it characteristic $1$-cycle} and a {\it reduced characteristic $1$-cycle} respectively. Note
that $Y_b$ is a disjoint union of the reduced characteristic $1$-cycles on
$Y_b$.

\medskip
{\bf Theorem 1.4} {\it  Under the same notations as in Theorem
1.3, let $Y_b$ be a reduction of a general singular fiber $M_b$.
Then the natural ${\mathbf C}^{n-1}$-action on $Y_b$ (cf. Proposition 2.2)
transitively acts on the set of characteristic $1$-cycles on $Y_b$. In particular, the characteristic $1$-cycles on $Y_b$
are all isomorphic. Moreover, modulo common divisors of $r_s$,
the characteristic $1$-cycle coincides with either

(1) one of singular fibers of a relatively minimal elliptic fibration
listed by Kodaira [Kd,
Theorem 6.2];

(2) $1$-cycle of Type $A_{\infty}$, i.e., $1$-cycle $\sum_{i \in
{\mathbf Z}} C_i$ consisting of infinitely many ${\mathbf P}_1$'s
such that $C_i \cap C_{i+1} = \{P_i\}$ (transversal and $P_i \not=
P_j$ if $i \not= j$), and such that $C_{i} \cap C_{j} = \emptyset$
if $\vert i - j \vert \ge 2$;

(3) $1$-cycle of Type $D_{\infty}$, i.e., $1$-cycle $C_{0} + C_{1}
+ \sum_{i \ge 2} 2C_i$ consisting of infinitely many ${\mathbf
P}_1$'s such that $C_{i} \cap C_{i+1} = \{P_i\}$ for each $i \ge
1$, $C_{0} \cap C_{2} = \{P_0\}$ (all are transversal and $P_i
\not= P_j$ if $i \not= j$) and such that $C_{i} \cap C_{j} =
\emptyset$ for other pairs $i \not= j$. }

\medskip
In other words, a general singular fiber is a disjoint union of a
family of one of Kodaira's singular fibers in [Kd, Theorem 6.2] or
an infinite cycle described in (2) or (3).  As Matsushita pointed
out to us, the case (2) actually occurs. His example will be given
in Proposition 4.13 in Section 4. Unfortunately, we do not have a
concrete example of the case (3), although we expect it exists.

\medskip
When $f$ is a projective morphism, Theorem 1.3 and some part of
Theorem 1.4  follow from the works of Matsushita ([M1], [M2]). In
particular, Matsushita gave a more or less complete classification
when $n=2$ and $f$ is projective in [M2], which gives more refined
information than ours, especially about the multiplicities and
monodromy. However, the view-point of the characteristic $1$-cycles
in Theorem 1.4 did not appear in his classification. We believe that
our view-point gives a new perspective even in the situation of [M2].

The argument used by Matsushita requires the projectivity of $f$,
because it depends on the classification theory of degeneration of
abelian varieties. No analog of this theory is known for
non-algebraic complex tori. Our proof of Theorems 1.3 and 1.4 uses
a completely different approach and uses the properties of
Lagrangian fibration more directly. A key tool of our approach is
twisted vector fields on the singular fibers. By examining the
Chern numbers of the leaves and the degree of the twisting, we can
control the degeneration of the fibers. This idea goes back to
Siu's work [Si] where he used certain twisted vector fields to
control the degeneration of complex structures. Another key
ingredient in the control of the singularity of a general singular
fiber is the theory of dualizing sheaves for singular curves, in
particular, the classical result of Rosenlicht's in [Se, Chapter
IV].

\medskip
In Section 2, we will introduce the basic geometric objects
associated to a proper holomorphic Lagrangian fibrations, called
the characteristic foliation of a vertical  hypersurface, and show
that Theorem 1.2 follows from Theorem 1.3. In Section 3, we study
the characteristic foliation arising from the determinantal
hypersurface more closely and prove Theorems 1.3. In Section 4, we
will prove Theorem 1.4,  using the characteristic vector fields
and the theory of dualizing sheaves for singular curves.

\medskip
While writing this paper, we learned about a preprint of J. Sawon
[Sa] where Question 1.1 was studied from a different view-point.
In particular, he found many interesting examples for which the
characteristic foliations are algebraic. But there is no overlap
with our result. Also after we finished a preliminary version of
this paper, we received Matsushita's preprint [M3] which gives a
classification of the general singular fibers when $f$ is a
projective morphism, refining his previous works [M1] and [M2].
His result is more refined than ours when $f$ is projective. His
approach generalizes that of [M2] and is completely different from
ours.

\medskip
{\bf Acknowledgement} We would like to express our thanks to
Daisuke Matsushita for pointing out an error in the first  version
of this paper and for informing us of a very impressive example in
Proposition 4.13.

\section{Characteristic foliation of a vertical hypersurface}

Let $M$ be a connected (not necessarily compact) complex manifold of
dimension $2n$. Assume that there exists a symplectic form $\omega$
on $M$. This means that $\omega$ is a $d$-closed holomorphic 2-form
 on $M$ which is non-degenerate at every point of $M$. The pair
$(M, \omega)$ is called  a {\it holomorphic symplectic manifold} of
dimension $2n$. The symplectic form $\omega$ defines a natural
${\mathcal O}_M$-linear isomorphism $\iota_{\omega}: T^*(M)
\longrightarrow T(M)$, which is the {\it inverse} of the contraction
map $T(M) \longrightarrow T^*(M)$ defined by $v \mapsto \omega(v,
*)$.

\medskip
 Let $ Y \subset M$ be a reduced (not necessarily irreducible) hypersurface in
$M$ and let $N^*_Y$
  be the conormal sheaf of $Y \subset M$. This is an
invertible sheaf on $Y$; in fact, $$ N^*_Y \cong {\mathcal
O}(-Y)|_Y\,\, .$$ Via the natural map $N^*_Y \longrightarrow T^*(M)|_Y$, we
get a canonical section of $T^*(M) \otimes {\mathcal O}(Y)|_Y$
uniquely defined up to a non-zero multiplicative function. Via the
isomorphism $\iota_{\omega}: T^*(M) \longrightarrow T(M)$, this gives
rise to a non-zero section $$\lambda_Y \in H^0( Y, T(M) \otimes
{\mathcal O}(Y))\,\, . $$ Locally, $\lambda_Y$ is defined as follows. Let
$h$, or more precisely $h = 0$, be a local defining equation of the hypersurface $Y$ in an
open subset $U \subset M$. Then on $U \cap Y$,
$$ \lambda_Y = \iota_{\omega}(dh)\,\, {\rm ,i.e.,}\,\,
\omega(\lambda_Y, *) = dh(*)\,\, .$$ From this local description,
it is clear that the leaves of the twisted vector field
$\lambda_Y$ is tangent to $Y$, defining a foliation on the smooth
part of $Y_{\rm reg}$ of $Y$. This foliation of $Y_{\rm reg}$ is exactly the
characteristic foliation  of $Y$ defined in the introduction. The
twisted vector field $\lambda_Y$ will be called the {\it
characteristic vector field} on $Y$.

\medskip
Now let $B$ be a complex manifold of dimension $n$ and let $f: M
\longrightarrow B$ be a holomorphic Lagrangian fibration, i.e., $f$
is a proper morphism with connected fibers such that the underlying
variety of each fiber of $f$ is a Lagrangian subvariety of $M$.
Recall that a reduced (not necessarily irreducible) subvariety $L$
of $(M, \omega)$ is called a Lagrangian subvariety if $L$ is of pure
dimension $n = \dim M/2$ and $\omega \vert_{L_{\rm reg}} = 0$.

Given a holomorphic function $g$ on $B$, the holomorphic
vector field $\iota_{\omega}(f^* dg)$ on $M$ is called the {\it
Hamiltonian vector field} associated to $g$. It is well known that
the Hamiltonian vector field is complete, because $f$ is proper,
and it is tangent to the fibers of $f$, in the sense that it maps the defining ideal of the reduction of each fiber to itself. For two functions $g_1$ and
$g_2$, their Hamiltonian vector fields $\iota_{\omega} (f^* dg_1)$
and $\iota_{\omega} (f^* dg_2)$ commute, i.e.,
$[\iota_{\omega} (f^* dg_1), \iota_{\omega} (f^* dg_2)] = 0$ under
the Lie bracket (see e.g. [Ar]).

We will say that a {\it reduced} hypersurface $Y \subset M$ is {\it
vertical} with respect to $f$ if the set-theoretic image $f(Y)$ is
an {\it irreducible} hypersurface in $B$.  The  requirement of the
irreducibility of  $f(Y)$ is somewhat artificial, but will be very
convenient when stating our results below. Note that  $Y$ itself is
not necessarily irreducible. We will denote the restriction of $f$
to $Y$ by $f|_Y \; : Y \longrightarrow f(Y).$ Here and hereafter, we
regard the image $f(Y)$ as an analytic subset with its reduced
structure. But we regard the preimage $f^{-1}(Z)$ of an analytic
subset $Z$ as an analytic subspace of $Y$ in the {\it
scheme-theoretic sense}.

\medskip
{\bf Proposition 2.1} {\it Let $Y \subset M$ be a vertical
hypersurface with respect to a Lagrangian fibration $f:M
\longrightarrow B$.  Then the characteristic foliation on $Y_{\rm reg}$ are
tangent to fibers of $f|_Y$.}

\medskip
{\it Proof}. The statement is  local on $Y$. So we may assume that
$Y$ is irreducible and replace $B$ by a coordinate neighborhood of a
smooth point of the hypersurface $f(Y)$. Let $y \in Y$ be a smooth
point of $Y$ such that $f(y)$ is a smooth point of $f(Y)$. Let $h
\in {\mathcal O}_{M, y}$ be a local defining equation of $Y$ at $y$
and $g \in {\mathcal O}_{B, f(y)}$ be a local defining equation of
$f(Y)$ at $f(y)$. Then $$f^*g = \xi \cdot h^k \mbox{ for some } \xi
\in {\mathcal O}_{M, y}^* \mbox{ and an integer } k >0\,\,.$$ The
characteristic vector field of $Y$ is defined locally by the 1-form
$dh$. On the other hand,
$$f^* d g = h^{k} \; d \xi + k \xi h^{k-1}  \; dh\,\,.$$ The
Hamiltonian vector field $\iota_{\omega}(f^* d g)$ is tangent to
the fibers of $f$ and  vanishes on $Y$ with multiplicity $k-1$.
Thus the local vector field
$$ \frac{1}{h^{k-1}} \iota_{\omega} (f^* d g) $$ defined in a
neighborhood of $y$ in $M$ is also tangent to the fibers of $f$.
Therefore, its restriction to $Y$ will be tangent to fibers of $f|_Y$
as well.
But the restriction is just
$$\iota_{\omega}(\frac{1}{h^{k-1}} f^* d g)|_Y = \iota_{\omega}
(k \xi \; dh)|_Y.$$ The right hand side is exactly the
characteristic vector field on $Y$ multiplied by the non-vanishing
function $k \xi$. It follows that the characteristic vector field
of $Y$ is tangent to the fibers of $f|_Y$. $\Box$

\medskip
{\bf Proposition 2.2} {\it Let $Y \subset M$ be a vertical
hypersurface with respect to a Lagrangian fibration $f: M
\longrightarrow B$. Then there exist $n-1$ commuting holomorphic
vector fields $v_1, \ldots, v_{n-1}$ in a neighborhood of a
general fiber $Y_b$ of $f|_Y$, which are linearly independent at
every point of $Y_b$. These vector fields generate an effective
holomorphic action of the commutative additive complex Lie group
${\mathbf C}^{n-1}$ on a general fiber of $f|_Y$ such that each
orbit has dimension $n-1$. In particular, the normalization
$\hat{Y}_b$ of a general fiber $Y_b$ of $f|_Y$ is smooth and the
connected component of the singular locus of $Y_b$ is a variety of
dimension $n-1$ homogeneous under the action of ${\mathbf
C}^{n-1}$, so that ${\rm Sing}\, (Y_b)$ is a disjoint union of
$(n-1)$-dimensional tori. }

\medskip {\it Proof}.
Choose a general point $b$ of $f(Y)$. Let $z_1, \ldots, z_{n-1}$ be
local holomorphic functions in a neighborhood $U$ of $b $ in $B$
whose restriction to $f(Y)$ give holomorphic coordinates in a
neighborhood of $b$ in $f(Y)$. Consider the holomorphic map $\Psi: U
\longrightarrow {\mathbf C}^{n-1}$ defined by $\Psi:= ( z_1, \ldots,
z_{n-1})$. By shrinking $U$, we may assume that $\Psi|_{f(Y) \cap U}
: f(Y) \cap U \longrightarrow \Psi (U)$ is biholomorphic. A general
fiber of $\Psi \circ f$ is smooth by Sard's theorem. From the
assumption that $b$ is a general point of $f(Y)$, we may assume that
the fiber of $\Psi \circ f$ over $\Psi (b)$ is smooth. In other
words, the 1-forms $f^* d z_1, \ldots, f^* d z_{n-1}$, defined in a
neighborhood of the fiber $f^{-1}(b)$ are linearly independent at
every point of $Y_b$. Via the isomorphism $\iota_{\omega} : T^*(M)
\longrightarrow T(M)$, these 1-forms define Hamiltonian vector
fields
$$v_i : = \iota_{\omega}(f^* d z_i)$$ in a neighborhood of $Y_b$.
Recall that these Hamiltonian vector fields are tangent to the
fibers of $f$. Moreover, they are linearly independent at every
point of $Y_b$ and they commute:
$$[v_i, v_j]=0, \mbox{ for each } 1 \leq i, j \leq n-1\,\, .$$
Since $f$ is proper, these Hamiltonian vector fields are complete.
Thus, they generate a holomorphic action of the commutative complex Lie
group ${\mathbf C}^{n-1}$ on $Y_b$. Each point of $Y_b$ has an
$(n-1)$-dimensional (not necessarily closed) orbit under this
action, as $v_1, \ldots, v_{n-1}$ are linearly independent at each point
of $Y_b$.
 This ${\mathbf C}^{n-1}$-action can be lifted
to the normalization $\hat{Y_b}$ of $Y_b$, because any local
analytic automorphism of a complex analytic variety  can be lifted
to its normalization. Let ${\rm Sing}\,(\hat{Y}_b)$ be the
singular locus of $\hat{Y_b}.$ By the normality, ${\rm
Sing}\,(\hat{Y}_b)$ has dimension $\leq n-2$. But each orbit of
the ${\mathbf C}^{n-1}$-action on $\hat{Y}_b$ has dimension $n-1$,
as the natural lifts of $v_1, \ldots, v_{n-1}$ under the natural
map $\hat{Y}_b \longrightarrow M$, say, $\hat{v}_1, \ldots,
\hat{v}_{n-1}$, are linearly independent at each point of
$\hat{Y}_b$. It follows that ${\rm Sing}\,(\hat{Y}_b)$ must be
empty. $\Box$

\medskip
{\bf Proposition 2.3} {\it Let $Y \subset M$ be a vertical
hypersurface with respect to a Lagrangian fibration $f: M
\longrightarrow B$. Let $Y_b$ be a general fiber of $f|_Y$ and
$\hat{Y_b}$ be its normalization, which is smooth by Proposition
2.2. Let $\hat{v}_i$ be the vector field on $\hat{Y_b}$ which is
the natural lift of the vector field $v_i$ defined in Proposition
2.2. Then there exist $n-1$ holomorphic 1-forms $\varphi_1,
\ldots, \varphi_{n-1}$ on $\hat{Y_b}$ with the following
properties.

(i)  $\varphi_i (\hat{v}_j) = \delta_{ij} $ for each $1 \leq i, j
\leq n-1,$

(ii) $\varphi_1, \ldots, \varphi_{n-1}$ are linearly independent
at every point of $\hat{Y_b}$, and

(iii) each $\varphi_i$ annihilates the lift of the characteristic
vector field $\lambda_Y$ to $\hat{Y_b}$.}

\medskip
{\it Proof}. Let $\nu: \hat{Y} \longrightarrow Y$ be the
normalization of $Y$ and let $\mu: \hat{Y} \longrightarrow f(Y)$ be
the composition $\mu = f \circ \nu$. Let $b \in f(Y)$ be a general
smooth point of $f(Y)$. Then by Proposition 2.2,  $\hat{Y}_b=
\mu^{-1}(b)$ is smooth and therefore $\mu$ is a smooth morphism
around $\hat{Y}_b$. In particular, $\hat{Y}$ is smooth in a
neighborhood of $\hat{Y}_b$. Let $\frac{\partial}{\partial z_i} \in
T_b(f(Y))$ be the coordinate vector with respect to the coordinates
$z_1, \ldots, z_{n-1}$ on $f(Y)$ used in the proof of Proposition
2.2. For each $i= 1, \ldots, n-1$, we will define a non-zero 1-form
$\varphi_i$ on $\hat{Y}_b$ as follows. For each $z \in \hat{Y}_b$,
pick an element $u_i \in T_z(\hat{Y})$ such that $\mu_*(u_i) =
\frac{\partial}{\partial z_i}$, where $\mu_*: T_z(\hat{Y})
\longrightarrow T_b(f(Y))$ is the differential of $\mu$, which is
surjective. Then one can define the holomorphic $1$-form $\varphi_i$
on $\hat{Y}_b$ by setting for each tangent vector $v \in
T_z(\hat{Y}_b)$,
$$\varphi_i (v) \; := \; \omega( \nu_*(u_i), \nu_*(v))\,\, ,$$ where
$\nu_*: T_z(\hat{Y}) \longrightarrow T_{\nu(z)}(M)$ is the
differential of $\nu: \hat{Y} \longrightarrow M$. This definition
does not depend on the choice of $u_i$ because $f$ is a Lagrangian
fibration. From this and the definitions of $v_i$ and
$\iota_{\omega}$, (i) is immediate. (ii) is then a direct
consequence of (i), as $\hat{v}_1, \ldots, \hat{v}_{n-1}$ are
linearly independent at each point of $\hat{Y}_b$.

To check (iii), it suffices to check it at a general point of $\hat{Y}_b$.
Choose a point $z\in \hat{Y}_b$ over the smooth
locus of $Y$ and  let $v \in T_z(\hat{Y}_b)$ be a local lift of
the characteristic vector field $\lambda_Y$. By the definition of
$\varphi_i$,
$$\varphi_i (v) \; = \; \omega( \nu_* (u_i), \lambda_Y) $$ which
is zero from the definition of the characteristic foliation. $\Box$

\medskip
{\bf Proposition 2.4} {\it Let $Y \subset M$ be a vertical
hypersurface with respect to a Lagrangian fibration $f: M
\longrightarrow B$. Let $Y_b$ be a general fiber of $f|_Y$ and
$\hat{Y}_b$ be its normalization, which is smooth by Proposition
2.2. Assume that $\hat{Y}_b$ is of class ${\cal C}$ so that Hodge
decomposition holds. Let $\hat{X}$ be an irreducible component of
$\hat{Y}_b$. Suppose that $ \dim H^0(\hat{X}, T^*(\hat{X})) \le
n-1$, where $2n = \dim M$. Then, the Albanese map $\alpha : \hat{X}
\longrightarrow {\rm
Alb}\, (\hat{X})$ is surjective, has connected fibers, and the
fibers of $\alpha$ correspond to the leaves of the characteristic
foliation on $Y$.}

\medskip
{\it Proof}. Recall that the holomorphic $1$-forms $\varphi_1, \ldots,
\varphi_{n-1}$ on $\hat{X}$ (in Proposition 2.3) are linearly independent
at {\it each} point of $\hat{X}$. Therefore, the dimension of
the image of the Albanese map is at least $n-1$. Thus $\alpha$ is surjective
 by our assumption. We can say more. Let $\beta : \hat{X} \longrightarrow A$ be the Stein factorization of $\alpha$. From
Proposition 2.3 (ii), the induced morphism $A \longrightarrow {\rm
Alb}\,(\hat{X})$ is \'etale and surjective. It follows that $A$ is
also an $(n-1)$-dimensional complex torus, and that
$H^{0}(\hat{X}, T^*(\hat{X})) = \beta^{*}(H^{0}(A, T^*(A)))$.
Thus, $\alpha$ coincides with $\beta$. Hence the general fibers of
Albanese map are connected algebraic curves on $\hat{X}$. They
correspond to the leaves of the characteristic foliation by
Proposition 2.3 (iii). $\Box$

\medskip
Note that by Proposition 2.4,  Theorem 1.2 follows from Theorem
1.3.

\section{Characteristic foliation of a discriminantal hypersurface}

Let $f:M \longrightarrow B$ be a holomorphic Lagrangian
fibration of class ${\cal C}$. Here $M$ and $B$ are complex
manifolds of dimension $2n$ and $n$. Consider the set of the
critical values of $f$, i.e.,
$$\Delta\,\, = \,\, \{ b \in B, f^{-1}(b) \mbox{ is singular} \}\,\, .$$
In this section, we study the characteristic foliation of the
special vertical hypersurface called the discriminantal
hypersurface more closely. First of all, we notice the following.

\medskip
{\bf Proposition 3.1}
\par
{\it (1) $M_s = f^{-1}(s)$ is an $n$-dimensional complex torus if
$s \not\in \Delta$.

(2) The critical set $\Delta$ is a hypersurface of $B$ unless $\Delta = \emptyset$.}

\medskip
The second statement is non-trivial. In fact, there is a flat morphism from
a $4$-dimensional manifold to a $2$-dimensional manifold whose critical set
consists of just one point [Ra].

\medskip
{\it Proof}. (1) is well-known as Liouville's theorem for a
Lagrangian fibration (e.g., [Ar]). The proof there is also valid in
the holomorphic symplectic category. Here we shall give another
proof with a similar flavor to the one of Proposition 2.4. By the
symplectic form $\omega$, we have
$$T^*(M_s) \cong N_{M_s} \cong {\mathcal O}_{M_s}^{\oplus n}\,\, .$$
Thus, $M_s$ has $n$ holomorphic $1$-forms which are linearly independent at
 each point of $M_s$. Since $M_s$ is of class $\mathcal C$, the Hodge
 decomposition theorem holds for $M_s$. Thus, as in the proof of Proposition 2.4,
 one can show that the Albanese map $ M_s \longrightarrow {\rm Alb}(M_s)$ is an isomorphism. This implies (1).

Let us show (2). It suffices to show that $f : M \longrightarrow B$
is a smooth morphism if $f$ is smooth over $B \setminus Z$ for some
analytic subset $Z$ of codimension $\ge 2$. Let $b \in Z$. Since the
statement is local at each $b$, we may (and will) assume that $B$ is
a germ of $b$ and shrink it whenever it is more convenient. Choose a
smooth point $P$ of $Y_b$. Then, there is a smooth $n$-dimensional
manifold $U \subset M$ such that $U$ meet $Y_b$ at $P$
transversally. The induced map $f \vert_U : U \longrightarrow B$ is
a finite map which is unramified over $B \setminus Z$. Since $Z$ is
of codimension $\ge 2$, by the purity of branched loci for a finite
map, $f \vert_U$ is unramified over $B$. Thus $U$ is a section of
$f$. In particular, $f^{-1}(b)$ is reduced at $P$. We denote $u_t =
(f \vert_U)^{-1}(t)$ ($t \in B$).

Let $z_1, \ldots , z_n$ be a local coordinate system of $B$ at $b$.
Using these coordinates, we obtain $n$ Hamiltonian vector fields,
$\iota_{\omega}(f^*dz_1), \ldots, \iota_{\omega}(f^*dz_n)$, which
are commutative, tangent to each fiber and linearly independent at
each point of $M$. Thus, as in the proof of Proposition 2.2, they
generate a holomorphic action of the commutative Lie group ${\mathbf
C}^n$ on $M$ over $B$ and define a holomorphic map
$$\rho : {\mathbf C}^n \times B \longrightarrow M\,\, ;\,\, (g, t) \mapsto g(u_t)\,\, .$$
For each $t \in B$, we denote the stabilizer subgroup
$\{g \in {\mathbf C}^n\, ;\, g(u_t) = u_t\,\}$ by $\Lambda_t$.
If $t \in B \setminus Z$, then $M_t = {\mathbf C}^{n}/\Lambda_t$,
as $M_t$ is $n$-dimensional complex torus by (1). Thus, we can naturally
regard $\Lambda_t$ as $H^1(M_t, \mathbf Z)$ (more precisely, the homology
group $H_1(M_t, \mathbf Z)$, but this does not matter).

Since $B \setminus Z$ is simply-connected, the local constant
systems $R^1(f \vert_{M \setminus f^{-1}(Z)})_* \mathbf Z$ is
actually constant; $R^1(f \vert_{M \setminus f^{-1}(Z)})_* {\mathbf
Z} \cong {\mathbf Z}^{2n}$ over $B \setminus Z$. The canonical basis
$e_1, \ldots, e_{2n}$ of the constant sheaf ${\mathbf Z}^{2n}$ gives
the generators of $H^{1}(M_t, \mathbf Z)$ and therefore the
generators of $\Lambda_t$ simultaneously for $t \in B \setminus Z$.
In this way, we obtain $2n$ ${\mathbf C}^{n}$-valued holomorphic
functions $f_1(t), \ldots , f_{2n}(t)$ on $B \setminus Z$ such that
$$\Lambda_{t} = \langle f_1(t)\, , \ldots\, , f_{2n}(t) \rangle$$ for each
$t \in B \setminus Z$.

Let $\Gamma := \pi_{1}(M \setminus f^{-1}(Z))$. Since ${\mathbf C}^n \times (B \setminus Z)$ is simply connected, this means
that $\Gamma \simeq {\mathbf Z}^{2n}$ and it acts as the holomorphic deck transformation on the universal cover ${\mathbf C}^n \times (B \setminus Z)$ over $B \setminus Z$. The action is the additive group action of the stabilizer subgroup $\Gamma_t$ on each ${\mathbf C}^{n} \times \{t\}$ ($t \in B \setminus Z$).
Since ${\mathbf C}^n \times Z$ is of codimension $\ge 2$ in ${\mathbf C}^n \times B$, the action $\Gamma$ uniquely extends to the action on
${\mathbf C}^n \times B$ over $B$. The induced action on
${\mathbf C}^{n} \times \{t\}$ ($t \in B$) is an additive group action.
({\it At the moment, we do not know if it is faithful and
discontinuous or not.})

Let $\tilde{M}$ be the universal covering space of $M$. As $f^{-1}(Z)$ is of codimension $\ge 2$, we have $\pi_{1}(M) \cong \Gamma$ and $\Gamma$ acts on the complex manifold $\tilde{M}$ as a holomorphic deck transformation. Since ${\mathbf C}^n \times (B \setminus Z)$ is simply connected, the natural inclusion $M \setminus f^{-1}(Z) \subset M$ lifts to a natural inclusion ${\mathbf C}^n \times (B \setminus Z) \subset \tilde{M}$, which is compatible with $\Gamma$. Let $\tilde{f} : \tilde{M} \longrightarrow B$ be
the composition of $f$ and the universal covering map. By the shape of the action of $\Gamma$, $\tilde{f}$ coincides with the second projection ${\mathbf C}^n \times (B \setminus Z) \longrightarrow B \setminus Z$ on ${\mathbf C}^n \times (B \setminus Z)$.

Since ${\mathbf C}^n \times B$ is also simply-connected, the
morphism $\rho : {\mathbf C}^n \times B \longrightarrow M$ lifts to
the morphism $\tilde{\rho} : {\mathbf C}^n \times B \longrightarrow
\tilde{M}$. $\tilde{\rho}$ is equivariant under the action of
$\Gamma$ and commutes with the morphisms to $B$, as so is on the
common open set ${\mathbf C}^n \times (B \setminus Z)$. Therefore
the additive action $\Gamma$ on ${\mathbf C}^n \times B$ is free and
discontinuous, as so is on $\tilde{M}$.

Thus, we obtain a smooth complex torus fibration $f' : ({\mathbf
C}^n \times B)/\Gamma \longrightarrow B$ and a bimeromorphic
morphism $\tau : ({\mathbf C}^n \times B)/\Gamma \longrightarrow M =
\tilde{M}/\Gamma$ over $B$. This morphism is also finite, as both
$({\mathbf C}^n \times B)/\Gamma$ and $M$ are proper over $B$ of
equidimensional fibers. Thus, $\tau$ is an isomorphism, as both
$({\mathbf C}^n \times B)/\Gamma$ and $M$ are normal (actually
smooth). Hence, $f$ is a smooth complex torus fibration as well.
$\Box$

\medskip
In what follows, we always assume that $\Delta \not= \emptyset$, and regard
$\Delta$ as a reduced hypersurface of $B$. We call $\Delta$ the {\it critical hypersurface} of $f$. We will fix one irreducible component $D$ of $\Delta$.
The scheme-theoretic preimage
$f^{-1}(D)$ is a divisor on $M$. We call $f^{-1}(D)$ the {\it discriminantal hypersurface} of $f$. Let $Y$ be the underlying reduced hypersurface of
$f^{-1}(D)$, i.e., $Y = f^{-1}(D)_{\rm red}$. Then $Y$ is a vertical hypersurface in the sense of Section 2. Let $ Y_b$ be a fiber of $f|_Y$ over a general point $b \in D$. Note that $Y_b$ is reduced (and of pure dimension $n$), as $b$ is general.

\medskip
{\bf Proposition 3.2} {\it $Y_b$ can not be a complex torus (See
Prop 3.5 for a stronger statement).}

\medskip
{\it Proof}. It is well-known that deformation of a Lagrangian
complex torus in a holomorphic symplectic manifold is
unobstructed, locally forming a fibration  over an $n$-dimensional
base space (e.g., [DM, Theorem 8.7]).  Thus if $Y_b$ is complex
torus, the scheme theoretic fiber $f^{-1}(b)$ must coincide with
(its reduction) $Y_b$, a contradiction to $b \in \Delta$. $\Box$

\medskip
Let $\nu: \hat{Y}_b \longrightarrow Y_b$ be the normalization of
$Y_b$. Then $\hat{Y}_b$ is a compact complex manifold by Proposition
2.2. Let $L$ be a line bundle on $M$. Let $w \in H^0(Y_b, T(M)
\otimes L)$ be a twisted vector field which is tangent to $Y_b$, in
the sense that its local expression as a vector field of $M$ around
each point $y \in Y_b$ maps the defining ideal of $Y_b$ to itself.
The local analytic automorphism of $M$ (near $y \in Y_b$)
corresponding to $w$ acts on $Y_b$ as a local analytic automorphism
of $Y_b$ (near $y$). Since any local analytic automorphism of a
complex analytic variety can be lifted to its normalization (and
$\hat{Y}_b$ is smooth by Proposition 2.2), $w$ can be lifted to a
twisted vector field
$$\hat{w} \in H^0(\hat{Y}_b, T(\hat{Y}_b) \otimes \nu^*L)\,\, .$$
Here, the existence of the smooth ambient manifold $M$ and the
fact that $L$ is a line bundle defined (not only on $Y_b$ but
also) on $M$ are important for the existence of $\hat{w}$ on
$\hat{Y}_b$.

\medskip
The next proposition is very crucial in the sequel.

\medskip
{\bf Proposition 3.3} {\it Let $L$ be a line bundle on $M$ and $w \in
H^0(Y_b, T(M) \otimes L)$. Assume that $w$ is tangent to $Y_b$
and that the lift $\hat{w}$
is annihilated by the 1-forms $\varphi_1, \ldots, \varphi_{n-1}$
of Proposition 2.3. Then $\hat{w}$ must vanish on the points of
$\hat{Y}_b$ lying over the singular locus of $Y_b$.}

\medskip
{\it Proof}. Since $w$ is locally defined in a neighborhood of
$Y_b$ and is tangent to $Y_b$,  $w$ must be tangent to the
singular locus of $Y_b$. By Proposition 2.2,  each component $E$
of the singular locus of $Y_b$ is a homogeneous variety under the
group action generated by commuting Hamiltonian vector fields
$v_1, \ldots, v_{n-1}$. Thus the  set $\nu^{-1}(E) $ is a smooth
hypersurface in $\hat{Y}_b$ and the lifted vector fields
$\hat{v}_1, \ldots, \hat{v}_{n-1}$ generate the tangent spaces of
$\nu^{-1}(E)$ at every point.  Suppose $\hat{w}$ does not vanish
at some point $z \in \nu^{-1}(E)$. Then
$$\hat{w}_z = c_1 \hat{v}_1 + \cdots + c_{n-1} \hat{v}_{n-1}$$ for
some complex numbers $ c_1, \ldots, c_{n-1}, $ one of which, say
$c_1$ is non-zero. This means that $\varphi_1(\hat{w}) $ is not
identically zero, and $\hat{w}$ is not annihilated by $\varphi_1$,
a contradiction. $\Box$

\medskip
Let $X$ be an irreducible component of $Y_b$. Let $k$ be the
multiplicity of $f^{-1}(b)$ at a general point of $X$. In other
words, the divisor $f^{-1}(D)$ has multiplicity $k$ along the
irreducible component $H$ of $Y$ such that $X \subset H$. Let
$g$ be a local defining equation of $f(Y)$ at $b$. Then $f^* d g$
vanishes on $H$ with multiplicity $k-1$. The vector field
$\iota_{\omega}(f^* d g)$ on $M$ vanishes on $H$ with multiplicity
$k-1$ as well. We can then consider the ${\mathcal O}(-(k-1) H)$-valued vector
field $\gamma$ defined locally by
$$\gamma := \frac{1}{h^{k-1}} \iota_{\omega}(f^* d g)$$ where $h$ is a local defining equation
of the reduced hypersurface $H$ in $M$. Then, globally, $$\gamma \in H^0(Y_b,
T(M) \otimes {\mathcal O}( -(k-1) H))\,\,.$$

\medskip
{\bf Proposition 3.4} {\it The twisted vector field $\gamma$
defined on $Y_b$ as above is proportional to the characteristic
vector field $\lambda_Y$ on the components of $Y_b$ where $\gamma$
is not identically zero.}

\medskip
{\it Proof}.  It is clear that $\gamma$ is proportional to $\lambda_Y$
on a component of $Y_b$ not lying on $H$. Let $X$ be a component of
$Y_b$, which lies on $H$. At a general
point $y \in X$, we have
$$f^*g = \xi \cdot h^k \mbox{ for some } \xi \in {\cal O}_{M,
y}^*\,\, , $$ as in the proof of Proposition 2.1. The local expression
of $\gamma$ is $ \frac{1}{h^{k-1}} \iota_{\omega} (f^* d g)$ near $y$.
Thus, $\gamma$ is proportional to $\lambda_Y$, as we have seen
it in the proof of
Proposition 2.1.  $\Box$

\medskip
{\bf Proposition 3.5} {\it Let $X$ be an irreducible component of
$Y_b$ and let $\hat{X}$ be its normalization. Then $ \dim
H^0(\hat{X}, T^*(\hat{X})) = n-1.$ In particular, the Albanese map
$\alpha: \hat{X} \longrightarrow {\rm Alb}\,(\hat{X})$ is
surjective and has connected fibers. Consequently, $X$ cannot be a
complex torus and  the characteristic foliation of $Y$ have
algebraic leaves by Proposition 2.4. }

\medskip
{\it Proof}. Let $\nu: \hat{X} \longrightarrow X$ be the normalization
map.  Since $\gamma \in H^0(Y_b, T(M) \otimes {\mathcal O}(-(k-1) H))$
and $\gamma$ is proportional to $\lambda_Y$ from Proposition 3.4,
we can apply Proposition 3.3 for $\gamma$, as the lift $\hat{\lambda}_Y$
of $\lambda_Y$ is annihilated by $1$-forms
$\varphi_{1}, \ldots , \varphi_{n-1}$ constructed in
Proposition 2.3. Thus the lifted vector field
$\hat{\gamma} \in H^0(\hat{X}, T(\hat{X}) \otimes \nu^*{\mathcal
O}(-(k-1)H)) $ on $\hat{X}$ vanishes on the points lying over the
singular locus of $Y_b$. We have also
the characteristic vector field $\lambda_H
\in H^0(X, T(M) \otimes {\mathcal O}(H))$ of the vertical hypersurface $H$,
which can be lifted to
$\hat{\lambda}_H \in H^0( \hat{X}, T(\hat{X}) \otimes \nu^* {\mathcal
O}(H)).$ We will use these two twisted vector fields
$\hat{\gamma}$ and $\hat{\lambda}_H$ on $\hat{X}$ to prove
Proposition 3.5.

Let $\varphi$ be a non-zero 1-form on $\hat{X}$. We need to show
that $\varphi$ is linearly dependent on $\varphi_1, \ldots,
\varphi_{n-1}$.

\medskip
{\bf Claim} {\it At every point $z \in \hat{X}$, $\varphi_z$ is
linearly dependent on $\varphi_1, \ldots, \varphi_{n-1}.$}

\medskip
{\it Proof of Claim}. Suppose to the contrary that at a general
point $z \in \hat{X}$, $\varphi_z$ is linearly independent from
$\varphi_1, \ldots, \varphi_{n-1}.$ Then both
$\varphi_z(\hat{\lambda}_H)$ and $\varphi_z(\hat{\gamma})$ are
non-zero. Thus we have non-zero sections
$$\varphi(\hat{\lambda}_H) \in H^0(\hat{X}, \nu^* {\mathcal O}(H)), \;
\varphi(\hat{\gamma}) \in H^0(\hat{X}, \nu^* {\mathcal O}(-(k-1)H))\,\, .$$
Then
$\varphi(\hat{\lambda}_H)^{\otimes (k-1)} \cdot
\varphi(\hat{\gamma})$ is a non-zero holomorphic function on $\hat{X}$.
But $\hat{\gamma}$ vanishes at a point of $\hat{X}$ over singular locus
of $Y_b$. This leads a contradiction {\it unless $X$ does not
contain any singular point of $Y_b$}.

From now, we assume furthermore that $X$ does not contain any
singular point of $Y_b$. Note then that $\hat{X} = X = Y_b$ and
$Y_b$ is smooth. If $k >1$, then both ${\mathcal O}(H)$ and
${\mathcal O}(-(k-1)H)$ have non-zero sections on $Y_b$. Thus
${\mathcal O}(H)$ is a trivial bundle on $Y_b$. This implies that
$\gamma \in H^0(Y_b, T(Y_b))$. Then $\gamma, v_1, \ldots, v_{n-1}$
are vector fields on $Y_b$ which are linearly independent at $z$.
Here $v_{1}, \ldots, v_{n-1}$ are vector fields in Proposition 2.2.
If a linear combination of these $n$ vector fields vanish at some
point of $\hat{X} = Y_b$, the holomorphic function $\varphi(\gamma)
= \varphi(\hat{\gamma})$ is non-constant, giving a contradiction.
Thus, the tangent bundle of $Y_b$ is trivial. In particular, there
are global $1$-forms $\eta_1, \ldots , \eta_n$ which are linearly
independent at every point of $Y_b$. Then the Albanese map of
$\alpha : Y_b \longrightarrow {\rm Alb}\, (Y_b)$ is \'etale and
surjective; It is in fact an isomorphism as in the proof of
Proposition 2.4. Thus, $Y_b$ is a complex torus, a contradiction to
Proposition 3.2. This completes the proof of Claim. $\Box$

\medskip
By this claim, at every point $z \in \hat{X}$, $\varphi_z$ is
linearly dependent on $\varphi_1, \ldots, \varphi_{n-1}.$ Since
$\varphi_1, \ldots, \varphi_{n-1}$ are linearly independent at
every point of $\hat{X}$, for each $y \in \hat{X}$, we can write
$$ \varphi_y = c_{1,y} \varphi_{1,y} + \cdots + c_{n-1, y} \varphi_{n-1, y}$$
for uniquely determined complex numbers $c_{1, y}, \ldots, c_{n-1,
y}$. Then $c_{i,y}$ is a holomorphic function on $\hat{X}$ and
must be independent of $y$.  It follows that $\varphi$ is linearly
dependent on $\varphi_1, \ldots, \varphi_{n-1}$ as 1-forms on
$\hat{X}$. $\Box$

\medskip
Now let us study the fibers of the Albanese map of $\hat{X}$.

\medskip
{\bf Proposition 3.6} {\it Let $C$ be a general fiber of the
Albanese map for $\hat{X}$. Then the ${\mathbf C}^{n-1}$-action on
$\hat{X}$ generated by $\hat{v}_1, \ldots, \hat{v}_{n-1}$ (cf.
Proposition 2.2) must act transitively on ${\rm Alb}\,(\hat{X})$ and
$\hat{X}$ is a holomorphic fiber bundle over ${\rm Alb}\,(\hat{X})$ with
$C$ as the fiber.}

\medskip
{\it Proof}. By the universality of the Albanese map, the group
${\mathbf C}^{n-1}$ of automorphisms of $\hat{X}$ acts
equivariantly on ${\rm Alb}\,(\hat{X})$. The action of ${\mathbf
C}^{n-1}$ on $\hat{X}$ is the one corresponding to the global
vector fields $\hat{v}_1, \ldots, \hat{v}_{n-1}$, which are
pointwisely dual to the basis $\varphi_{1}, \ldots, \varphi_{n-1}$
of $H^{0}(\hat{X}, T^*(\hat{X}))$ (Proposition 2.3). Thus, from
the construction of Albanese variety, the equivariant action of
${\mathbf C}^{n-1}$ is the natural translation on the
$(n-1)$-dimensional complex torus ${\rm Alb}\, (\hat{X})$. This
implies the result. $\Box$

\medskip
The next proposition describes the structure of a smooth fiber
$Y_b$.  This proposition says that a smooth fiber is nothing but a
hyperelliptic surface when $\dim\, M = 4$ (cf. [M2], Table 4 Type
$I_{0}$).

\medskip
{\bf Proposition 3.7} {\it Suppose $Y_b$ is smooth. Then the
Albanese map of $Y_b (= \hat{X})$ is an elliptic fiber bundle.
Moreover, there exists an action of an $(n-1)$-dimensional torus
on $Y_b$ which acts transitively on ${\rm Alb}\,(Y_b)$  rendering
$Y_b$ the structure of a Seifert fibration over ${\mathbf P}_1$
(see e.g. [Ho] for the definition).}

\medskip
{\it Proof}. If $Y_b$ is smooth, then it is irreducible and $\hat{X}
= X= Y_b$. We have $Y=H$, $k>1$ and ${\cal O}( k H)$ is a trivial
bundle on $Y_b$. Thus $C \cdot H =0$ for the fiber $C$ of the
Albanese map. On the other hand, the characteristic vector field
$\lambda_H$ gives rise to a non-zero element of $H^{0}(C, T(C)
\otimes {\mathcal O}(H))$. Thus ${\rm deg} (T(C) \otimes {\mathcal
O}(H)) \ge 0$. Therefore ${\rm deg}\, T(C) \ge 0$ by $C \cdot H =
0$. It follows that $C$ is either a smooth elliptic curve or
${\mathbf P}_1$. Suppose $C = {\mathbf P}_1$. The hypersurface $Y$
on $M$ is smooth in a neighborhood of $C$ and the normal bundle of
$C$ in $Y$ is trivial. Moreover, the normal bundle of $Y$ in $M$ is
just ${\mathcal O}(H)$. This implies that the normal bundle of $C $
in $ M$ is trivial by $C \cdot H = 0$. But then $K_M \cdot C = {\rm
deg}\, K_{C} = -2$ by the adjunction formula. This is a
contradiction to the fact that $K_M$ is trivial by the symplectic
form $\omega$. Thus $C$ is an elliptic curve. This means that $Y_b$
contains no rational curve, as the Albanese map is a fiber bundle by
Proposition 3.6. Applying [Fu,Theorem 5.5], we see that ${\rm
Aut}_o\,(Y_b)$ is an $(n-1)$-dimensional torus  isogenous to ${\rm
Alb}\,(Y_b)$. From the general result about torus action in [Ho], we
have a Seifert fibration $\sigma: Y_b \longrightarrow C'$ over a
smooth curve $C'$ such that fibers of $\sigma$ are ${\rm
Aut}_o\,(Y_b)$-orbits. If $C'$ has a non-trivial 1-form $\psi$, then
$\sigma^* \psi$ is a non-trivial 1-form on $Y_b$ which is
non-trivial on $C$ as well. This is a contradiction to the fact that
$C$ is a fiber of the Albanese map for $Y_b$. Thus $C' \cong
{\mathbf P}_1$. $\Box$

\medskip
{\bf Proposition 3.8} {\it Suppose $Y_b$ is singular. Then for each
component $X$ of $Y_b$, its normalization $\hat{X}$ is a ${\mathbf
P}_1$-bundle over the complex torus ${\rm Alb}\,(\hat{X})$. In
particular, each component of $Y_b$ is covered by rational curves
and each rational curve on $Y_b$ is the image of the Albanese fiber
of $\hat{X}$ for some component $X$.}

\medskip
{\it Proof}. We will use the notation in the proof of Proposition
3.5. The fiber $C$ of the Albanese map of $\hat{X}$ is a leaf of
the two twisted vector fields $$\hat{\lambda}_H \in H^0( \hat{X},
T(\hat{X}) \otimes \nu^* {\mathcal O}(H)), \;\; \hat{\gamma} \in
H^0(\hat{X}, T(\hat{X}) \otimes \nu^*{\mathcal O}(-(k-1)H))\,\, . $$

Moreover, by Propositions 3.3, both $\hat{\lambda}_H$ and
$\hat{\gamma}$ vanish on the  divisor $E$ lying over ${\rm Sing}\,
Y_b$. The divisor $E$ certainly meets $C$ properly by Propositions
2.2 and 3.6. Therefore, both twisted vector fields
$$0 \not= \hat{\lambda}_H \vert_{C} \in H^0(C,
T(C) \otimes \nu^* {\mathcal O}(H)), \;\; 0 \not= \hat{\gamma}\vert_{C} \in
H^0(C, T(C) \otimes \nu^*{\mathcal O}(-(k-1)H))$$ have non-empty
zeros. Hence
$${\rm deg}\, T(C) + \nu^*H \cdot C > 0, \, \;\; {\rm deg}\, T(C) -(k-1) \nu^*H \cdot C > 0\,\, .$$
Therefore, ${\rm deg}\, T(C) > 0$. Hence $C = {\mathbf P}_1$. A
rational curve on $Y_b$  cannot lie on ${\rm Sing}(Y_b)$ by
Proposition 2.2. Thus it must come from a rational curve on the
normalization, that is, it must be the image of the Albanese fiber
of $\hat{X}$ for some component $X$. $\Box$

\medskip
Propositions 3.5, 3.7 and 3.8 complete the proof of Theorem 1.3 and
consequently the proof of Theorem 1.2 via Proposition 2.4.

\section{Structure of general singular fibers}

In this section, we shall give a Kodaira-type classification of  a
general singular fiber of a holomorphic Lagrangian fibration $f :
M \longrightarrow B$ of class $\mathcal C$. Theorem 1.4 will
follow from Propositions 3.7, 3.8, 4.10, 4.11 and 4.12.

In what follows, {\it unless stated otherwise}, we use the same
notations as in Section 3. For instance, $2n = \dim\, M = \dim\,
B$; $b$ is a general point of an irreducible component $D$ of the
critical hypersurface $\Delta \subset B$; $Y$ is the underlying
reduced hypersurface of $f^{-1}(D)$; $Y_b$ is the fiber
$(f\vert_Y)^{-1}(b)$; $X$ is an irreducible component of $Y_b$;
$\nu : \hat{X} \longrightarrow X$ is the normalization of $X$;
$\alpha : \hat{X} \longrightarrow {\rm Alb}\,(\hat{X})$ is the
Albanese map (recall that $\hat{X}$ is smooth by Proposition 2.2
and that $\alpha$ is surjective with connected fibers by
Proposition 3.5); $H$ is the irreducible component of $f^{-1}(D)$
containing $X$, and $k$ is the multiplicity of $H$ in $f^{-1}(D)$.
Besides these notations, we denote by $C$ a fiber of the Albanese
map $\alpha : \hat{X} \longrightarrow {\rm Alb}\,(\hat{X})$. Note
that $\alpha$ is a holomorphic fiber bundle with typical fiber $C$
(Proposition 3.6) and that $C$ is a smooth elliptic curve when
$Y_b$ is smooth and $C = {\mathbf P}_1$ when $Y_b$ is singular
(Propositions 3.7 and 3.8).

Since Theorem 1.4 is local at $b \in D \subset B$ in the classical
topology, we will freely shrink $B$ around $b$, whenever it is
more convenient. For example, we may (and will) assume that $X = H
\cap Y_b$ for an irreducible component $X$ of $Y_b$.

\medskip
In Propositions 4.1-4.7 below, we study the structure of a germ of
$Y_b$ at its singular point. In Proposition 4.1, we consider the
case where the germ is not irreducible. We consider the case where
the germ is irreducible in Lemma 4.2-Proposition 4.7.

To state Proposition 4.1, it is convenient to define the following
notions. Let $x$ be a singular point of $Y_b$ and ${\mathcal X}_i,
i=0, \ldots, \ell$ be the irreducible components of the germ of
$Y_b$ at $x$. By Proposition 3.8, there exists a unique germ
${\mathcal C}_i \subset {\mathcal X}_i$ of a rational curve on
$Y_b$ through $x$ for each $i$. Suppose that for two components
${\mathcal X}_i$ and ${\mathcal X}_j, i \neq j,$  the
scheme-theoretic intersection ${\mathcal X}_i \cap {\mathcal
X}_j$, which is of dimension $n-1$ by Proposition 2.2,  defines a
Cartier divisor ${\mathcal D}_i$ on ${\mathcal X}_i$ and a Cartier
divisor ${\mathcal D}_j$ on ${\mathcal X}_j$. For example, this is
the case if both ${\mathcal X}_i$ and ${\mathcal X}_j$ are smooth.
We will say that ${\mathcal X}_i$ and ${\mathcal X}_j$ {\it
intersect transversally} if the local intersection numbers at $x$
are
$${\mathcal D}_i \cdot {\mathcal C}_i = {\mathcal D}_j \cdot
{\mathcal C}_j =1.$$ We will say that the two components {intersect
with multiplicity 2} if $${\mathcal D}_i \cdot {\mathcal C}_i =
{\mathcal D}_j \cdot {\mathcal C}_j =2.$$

\medskip
{\bf Proposition 4.1} {\it If the germ of $Y_b$ at a point is not
irreducible, then one of the following holds.

(Case 1) The germ has three irreducible components. Each component is smooth
and each pair intersects transversally.

(Case 2) The germ has two irreducible components. Each component is smooth
and the two components intersect with multiplicity 2.

(Case 3) The germ has two irreducible components. Each component is smooth
and the two components intersect transversally. }

\medskip
{\it Proof}. Let $x \in Y_b$ be such a point and $h$ be the local
defining equation of $D$ at $b$. Let
$$f^*h = h_0^{a_0} \cdots
h_{\ell}^{a_{\ell}}\,\, ,\,\, 0< a_0 \leq \cdots \leq a_{\ell}$$ be
a unique factorization into irreducible factors modulo units in
${\mathcal O}_{M, x}$. Put ${\mathcal H}_{i} = {\rm div}\, h_{i}$.
These ${\mathcal H}_{i}$ form local analytic irreducible components
of $Y = f^{-1}(D)_{\rm red}$ in a neighborhood $U$ of $x$ in $M$. We
choose as $H$ the global irreducible component of $Y$ such that
${\mathcal H}_{0} \subset H$. We put ${\mathcal X}_i = {\mathcal
H}_i \cap Y_b$. Let $X$ be a global irreducible component of $Y_b$
such that $$ x \in {\mathcal X}_0 \subset X \subset H.$$
 We also choose an
Albanese fiber $C$ of $\hat{X}$ such that $x \in \nu(C)$. Since
$f^{-1}(D)$ is a principal divisor on $Y_b$, it is trivial on
$\hat{X}$ as a line bundle. Thus, so is it on $C$ and $C \cdot
\nu^*(f^{-1}(D)) = 0$. Recall that $Y = f^{-1}(D)_{\rm red}$. Then
$C \cdot \nu^*Y \le0$, as $\nu(C) \subset X$ and $a_0 \leq \cdots
\leq a_{\ell}$. Since $C = {\mathbf P}_1$, one can then regard the
restriction $\hat{\lambda}_Y \vert_C$ of the lift of the
characteristic vector field $\lambda_Y$ as a non-zero vector field
on $C$ (and we shall do so). Again, since $C = {\mathbf P}_1$, the
vector field has at most two zero over $x$ counted with
multiplicities.  $C$ meets each divisor of $\hat{X}$ lying over
${\rm Sing}\, {\mathcal X}_0$ (if it is not empty) and also meets
the divisors lying over ${\mathcal X}_i \cap {\mathcal X}_0$ ($1
\le i \le \ell$). There, $\hat{\lambda}_Y \vert_C$ is zero with
multiplicities counted by the local expression $\iota_{\omega}(h_1
\ldots h_{\ell} dh_0)$. Hence $\ell \le 2$ and the Case (1) occurs
when $\ell=2$.

Consider the case where $\ell = 1$. If ${\mathcal X}_0$ is singular
at $x$, then $\nu^* h_1$ has two zeros counted with multiplicities
on $C$ and $\iota_{\omega}(dh_0)$ has additional zeros on $C$, a
contradiction. Thus ${\mathcal X}_0$ is smooth and ${\mathcal X}_0$
intersects the other component with multiplicities at most $2$ at
$x$. If $a_0 = a_1$, then, by changing the role of $h_0$ and $h_1$,
we see that the other local irreducible component is smooth as well.

It remains to consider the case where $\ell = 1$ and $a_0 < a_1$.
We use the twisted vector field $\lambda_H$ instead of $\lambda_{Y}$
above. By $a_0 \not= a_1$, the global irreducible component $H'$ of
$Y$ such that ${\mathcal H}_1 \subset H'$ is different from $H$.
Therefore, the global irreducible decomposition of $f^{-1}(D)$ is of
the form $a_0 H + a_1 H' + \cdots$. Since $C \cdot \nu^* H' > 0$ and
$0 < a_0 < a_1$, one has then
$$0 = C \cdot f^{-1}(D) \ge a_0 (C\cdot \nu^* H) + a_1 (C\cdot \nu^* H')
> a_0(C\cdot \nu^* H + C\cdot \nu^* H')\,\, .$$
Thus, if $C \cdot \nu^* H' \ge 2$, then $C \cdot \nu^* H < -3$. However,
this contradicts the fact that
$$0 \not= \hat{\lambda}_H \vert_C \in H^0(C, T(C) \otimes \nu^* {\mathcal O}(H))\,\, .$$
Thus $C \cdot \nu^* H' = 1$. In particular, the two local irreducible
components are smooth and meet transversally at $x$. This completes the proof.
$\Box$

\medskip
To study the irreducible germ of $Y_b$, we need some preliminary
results regarding dualizing sheaves. Let $V'$ be a complex variety
with the property that its normalization is smooth and the
dualizing sheaf $\omega_{V'}$ is invertible. Denoting by $\nu: V
\longrightarrow V'$ the normalization map, we have a natural
injective sheaf map $\omega_V \longrightarrow \nu^* \omega_{V'}$.
Since this is an inclusion of invertible sheaves on the smooth
variety $V$, we can naturally identify $\omega_V = \nu^*
\omega_{V'} \otimes {\mathcal O}_V( -E)$ for some effective
divisor $E$ on $V$. In Lemma 4.2 and Proposition 4.5 below, we
will study this effective divisor $E$ in some cases.

\medskip
{\bf Lemma 4.2} {\it Let $R'$ be an irreducible germ of an analytic
curve with an isolated singularity $Q$ and let $\nu: R \rightarrow
R'$ be the normalization morphism. Denote by $\omega_R$ and
$\omega_{R'}$ the dualizing sheaves of the curves. Assume that
$\nu^{-1}(Q) = P \in R$ (as a set) and that $\omega_{R'}$ is locally
free. Then  $$\omega_R =  \nu^* \omega_{R'} \otimes {\mathcal
O}_{R}(- 2\delta P)$$ where $\delta$ is the codimension of ${\cal
O}_{R', Q}$ in ${\cal O}_{R, P}$.}

\medskip
{\it Proof}. Since $\omega_{R'}$ is locally free and $\nu^{-1}(Q) =
P$, it follows from [Se, Page 72] that
$$\nu^* \omega_{R'} = {\mathcal O}_{R} \cdot
\frac{\eta}{t^{2\delta}}\,\, ,\,\, \omega_{R} = {\mathcal O}_{R}
\cdot \eta\,\, .$$ Here $\eta$ is a local generator of $\omega_R$
and $t$ is a generator of the maximal ideal ${\mathbf m}_{R, P}$.
Thus,
$$\nu^* \omega_{R'} = \omega_{R}(2\delta P) = \omega_{R} \otimes {\mathcal O}_{R}(2\delta P)\,\, .$$
Since $\omega_{R}$ is invertible, this gives the desired equality.  $\Box$

\medskip
{\bf Lemma 4.3} {\it Let $V'$ be an irreducible local complete
intersection variety in a smooth manifold $\mathcal V$ and let
$\rho: V' \rightarrow Z$ be a surjective morphism to an irreducible
variety $Z$. Then for a general fiber $R'$ of $\rho$, the dualizing
sheaf $\omega_{R'}$ is the restriction of the dualizing sheaf
$\omega_{V'}$ to $R'$. }

\medskip

{\it Proof}. Put $\ell = \dim\, Z$. Let $R':= \rho^{-1}(P)$ be a
fiber over a general smooth point $P$ of $Z$. Then $R'$ is
reduced. The ideal sheaf $I_{V', R'}$ of $R'$ in $V'$ is generated
by $\rho^* t_1. \ldots , \rho^*t_{\ell}$, where $t_1, \ldots,
t_{\ell}$ is a generator of the maximal ideal ${\mathbf m}_{Z,
P}$. The functions $\rho^* t_1. \ldots , \rho^*t_l$ form a regular
sequence in ${\mathcal O}_{V'}$. Thus $R'$ is a local complete
intersection in $\mathcal V$ as so is $V'$. Let $a$ and $b$ be the
codimensions of $V'$ and $R'$ in $\mathcal V$. Then we have
natural isomorphisms
$$\omega_{V'} \cong \omega_{\mathcal V}
\otimes \wedge^{a} (I_{{\mathcal V}, V'}/{I}_{{\mathcal V},
V'}^{2})^{*}\,\, ,\,\,  \omega_{R'} \cong \omega_{\mathcal V}
\otimes \wedge^{b} (I_{{\mathcal V}, R'}/{I}_{{\mathcal V},
R'}^{2})^{*}\,\, ,$$ as in [GPR, Chapter II, Section 5]. Since $R'$
is reduced, the smooth locus of $R'$ is a Zariski open dense subset.
Thus, the natural morphism
$$I_{{\mathcal V}, V'}/{I}_{{\mathcal V}, V'}^{2} \vert_{R'} \longrightarrow
I_{{\mathcal V}, R'}/{I}_{{\mathcal V}, R'}^{2}$$ is generically
injective. Indeed,  this is nothing but the natural injective
morphism of the usual conormal bundles, when we restrict it to the
smooth locus of $R'$. Since $I_{{\mathcal V}, V'}/{I}_{{\mathcal V},
V'}^{2} \vert_{R'}$ is locally free, the generically injective
morphism above has to be injective. Thus, the following natural
complex is exact:
$$0 \longrightarrow I_{{\mathcal V}, V'}/{I}_{{\mathcal V}, V'}^{2} \vert_{R'} \longrightarrow I_{{\mathcal V}, R'}/
{I}_{{\mathcal V}, R'}^{2} \longrightarrow I_{V', R'}/ {I}_{V',
R'}^{2} \cong {\mathcal O}_{R'}^{\oplus (\dim V' - \dim Z)}
\longrightarrow 0\,\, .$$ Here the last isomorphism follows from the
fact that $R'$ is a fiber over a smooth point $P$. Taking the dual
of this sequence and then the wedge product, we have a natural
isomorphism
$$\wedge^{a} (I_{{\mathcal V}, V'}/{I}_{{\mathcal V}, V'}^{2})^{*}  \vert_{R'} \cong \wedge^{b} (I_{{\mathcal V}, R'}/{I}_{{\mathcal V}, R'}^{2})^{*}\,\, .$$
This implies the result. $\Box$

\medskip
{\bf Proposition 4.4} {\it

(1)  Each irreducible component of $Y_b$, as well as $Y_b$ itself,
is a local complete intersection variety. Consequently, it has
invertible dualizing sheaf.

(2) Let $x$ be a point of $Y_b$. Then there exists a germ of an
analytic curve $R'$ in ${\bf C}^2$ and a biholomorphic map from
the germ of $Y_b$ at $x$ to the product $R' \times {\mathcal M}$
where ${\mathcal M}$ denotes the germ of an $(n-1)$-dimensional
complex manifold. Moreover the orbits of ${\bf C}^{n-1}$-action on
$Y_b$ correspond to the ${\mathcal M}$-factors.

(3) Let $x$ be a singular point of $Y_b$. Then there exist an open
neighborhood $U$ of $x$ in $M$ and a (not necessarily integrable)
distribution ${\mathcal D}$ of rank 2 on $U$ such that it is normal to
${\rm Sing}(Y_b)$ at $x$ and such that the germ of any rational curve
$C'$ of $Y_b$ through $x$ is tangent to ${\mathcal D}$ at smooth
points of $C'$. }

\medskip
{\it Proof}. Let $z_1 , \ldots , z_{n-1}$ be the same as in the
proof of Proposition 2.2. From the non-vanishing of $f^* d z_1,
\ldots, f^* d z_{n-1}$ in a neighborhood of $Y_b$, we see that $Z_b
:= f^{-1}(z_1= b_1, \ldots, z_{n-1}= b_{n-1})$ is smooth in a
neighborhood of $Y_b$ for general $b = (b_1, \ldots, b_{n-1}) \in D$
and $Y_b$ is just a reduced hypersurface in $Z_b$. This implies (1).

The commuting vector fields $v_1, \ldots, v_{n-1}$ of Proposition
2.2 are defined in a neighborhood of $Y_b$.  In a neighborhood $U$
of $x$ in $Z_b$, we can define coordinates $w_1, \ldots, w_{n+1}$
such that $v_i = \frac{\partial}{\partial w_i}$ for $1 \leq i \leq
n-1$. Since the hypersurface $Y_b$ in $Z_b$ is invariant under $v_1,
\ldots, v_{n-1}$, we can choose a local defining equation of $Y_b$
on $U$ as a holomorphic function depending only on the two variables
$w_n$ and $w_{n+1}$. This implies (2).

To see (3), it suffices to construct such a distribution on the
germ of $Z_b$ at $x$. We use the construction in the proof of
Proposition 2.3. Let $\frac{\partial}{\partial z_i}, 1 \leq i \leq
n-1,$ be the vector fields on $B$ transversal to $D$ near $b$ and
let $u_i$ be a vector field defined in a neighborhood of $x$ in
$M$ satisfying $f_*(u_i) = \frac{\partial}{\partial z_i}$ for each
$1 \leq i \leq n-1$. Let $\phi_i$ be the 1-form in this
neighborhood defined by $u_i = \iota_{\omega}(\phi_i)$. Note that
our $\varphi_i$ in Proposition 2.3 is the pull-back of $\phi_i$ to
the normalization of $Y_b$. Unlike $\varphi_i$,  $\phi_i$ is
defined only in a neighborhood of $x$ because its definition
depends on the choice of $u_i$.  These $n-1$ differential forms
restricted to the $(n+1)$-dimensional manifold $Z_b$ define a
distribution ${\cal D}$ of rank $2$ on $Z_b$. From the obvious
relation $\phi_i (v_j) = \delta_{ij}$, this distribution is normal
to the ${\bf C}^{n-1}$-orbits and consequently normal to ${\rm
Sing}(Y_b)$. Since
 the rational curve $C'$ is the fiber of the Albanese map defined by $\varphi_1,
 \ldots, \varphi_{n-1}$, $C'$ is tangent to ${\mathcal D}$ at its
 smooth point.
 $\Box$

\medskip
{\bf Proposition 4.5} {\it Let $X$ be a component of $Y_b$. Assume
that
 there exists a codimension-1 irreducible component $S'$ of the
singular locus of $X$ such that the germ of $X$ at a general point
of $S'$ is irreducible so that the support of $\nu^{-1}(S')$ is an
irreducible reduced divisor $S$ on $\hat{X}$ and $\nu|_S$ is
bimeromorphic. Then
$$ \omega_{\hat{X}} =
\nu^* \omega_{X} \otimes  {\mathcal O}(- 2 \delta S) \otimes
{\mathcal O}(-E')$$ for some positive integer $\delta$ and an
effective divisor $E'$ on $\hat{X}$. }

\medskip
{\it Proof}. Since the statement is local in a neighborhood of a
general point of $S'$, it suffices to prove $$ \omega_{\hat{X}} =
\nu^* \omega_{X} \otimes  {\mathcal O}(- 2 \delta S)$$
 for the germ of $X$
at a general point of $S'$.  By the description of the germ in
Proposition 4.4, this follows immediately from Lemma 4.2 and Lemma
4.3.  $\Box$

\medskip
{\bf Proposition 4.6} {\it Let $X$ be a component of $Y_b$ and let
$C$ be a fiber of the Albanese map of $\hat{X}$. Then

 (1) $C \cdot (\nu \circ f)^{-1}(D) = 0$ and $C \cdot \nu^* \omega_X = C \cdot
\nu^* H \le 0$.

(2) Assume that $Y_b$ is not irreducible. Then $C \cdot \nu^*
\omega_X = C \cdot \nu^* H \le -1$.}

\medskip
{\it Proof}. By $\omega_{M} \simeq {\mathcal O}_{M}$, we have
${\mathcal O}(H) \vert_H \simeq \omega_{H}$ by the adjunction
formula. Since we have assumed, by shrinking $B$, that $X$ is a
fiber of $f\vert_H : H \longrightarrow f(H)$, we have $\omega_{H}
\vert_{X} = \omega_{X}$ as in the proof of Lemma 4.3. Thus,
${\mathcal O}(H) \vert_{X} \simeq \omega_{X}$. This implies the
second equality in (1). We have $C \cdot \nu^*(f^{-1}(D)) = 0$ as
$f^{-1}(D)$ is a principal divisor on $Y_b$. We have also that $C
\cdot \nu^* H' \ge 0$ for each irreducible component $H'$ of $Y =
f^{-1}(D)_{\rm red}$ such that $H' \not= H$. Hence $C \cdot \nu^* H
\le 0$. This implies (1).

If $Y_b$ is not irreducible, then there is at least one irreducible component
$H'$ of $Y = f^{-1}(D)_{\rm red}$ such that $H' \not= H$.
Moreover, $C \cdot \nu^* H' > 0$ at least one of such $H'$, as $Y_b$ is connected. Therefore
$C \cdot \nu^* H \le -1$. $\Box$

\medskip
{\bf Proposition 4.7} {\it Suppose a component $X$ of $Y_b$ is
locally irreducible but singular at a point $x$. Then the following
holds.

(1) $Y_b$ is irreducible and the singular locus of $X$ is
irreducible.

(2) The germ of $X$ at $x$ is  biholomorphic to the product of the
singularity of the rational cuspidal cubic curve in the plane and a
complex manifold.

(3) Each rational curve on $Y_b$ is isomorphic to the rational
cuspidal cubic curve in the plane. }

\medskip
{\it Proof}. Let $C$ be a fiber of the Albanese map of $\hat{X}$
such that $x \in \nu(C)$. In the notation of Proposition 4.5
$$-2= \omega_{\hat{X}} \cdot C =
\nu^* \omega_{X} \cdot C - 2 \delta S \cdot C -E'\cdot C\, ,$$
which implies that $\nu^* \omega_{X} \cdot C = \nu^* {\mathcal
O}(H) \cdot C \geq 0.$  By Proposition 4.6 (1), we conclude that
$$\nu^* {\mathcal O}(H) \cdot C = 0\, ,\; \delta =1\, , \; S
\cdot C = 1\, , \; E' \cdot C =0\, . $$ This implies that $Y_b$ is
irreducible by Proposition 4.6 (2). Also, it implies that $\nu(S)$
is the only component of ${\rm Sing} (X)$ where $X$ is locally
irreducible, because otherwise $E' \cdot C >0$. In particular, if
there exists another component of ${\rm Sing}(X)$, there are two
distinct points $x_1, x_2$ on $C$ disjoint from $S$ such that
$\nu(x_1), \nu(x_2) \in {\rm Sing}(X)$. From $\nu^* {\mathcal
O}(H) \cdot C \geq 0$, $\nu^*{\cal O}(H) \vert_C$ is a trivial
line bundle on $C = {\mathbf P}_1$. The twisted vector field
$\hat{\lambda}_H \vert_C \not= 0$ is then a vector field on $C$
vanishing at the points $x_1, x_2$ and $\nu^{-1}(x)$, a
contradiction. Thus ${\rm Sing}(X)$ is irreducible. This completes
the proof of (1).

By Proposition 4.4 (2), we already know that the germ of $X$ at $x$
is of the form $R' \times {\mathcal M}$. Let us use the notation of
Lemma 4.2. From [Se, p. 59, equation (1)], we have the inclusion
$${\mathbf C} + c_{x} \subset {\mathcal O}_{R', x} \subset
{\mathcal O}_{R, \nu^{-1}(x)}\, ,$$ where $c_x$ is the conductor
of $R'$ at $x$. Since  $\delta =1,$ $c_x$ is the square of the
maximal ideal of $R$ at $\nu^{-1}(x)$ by [Se, p.71, Section 11].
Then $ {\mathbf C} + c_{x} = {\mathcal O}_{R', x}$ and $R'$ is the
germ of the rational cuspidal cubic plane curve at the cusp. This
completes the proof of (2).

For (3), it suffices to show that the germ of $C' = \nu(C)$ at $x$
is isomorphic to that of $R'$. Since $C$ is transversal to $S$, we
see that under the projection $R' \times {\mathcal M}
\longrightarrow R'$,  the germ of $C'$ is projected to $R'$
bijectively. From the property of the cusp of $R'$, either the germ
of $C'$ is isomorphic to $R'$ or $C'$ is non-singular. In the latter
case, $C'$ corresponds to a smooth curve on $R' \times {\mathcal
M}$, and must be tangent to ${\rm Sing}(Y_b)$ at $x$. However, from
Proposition 4.4 (3), if $C'$ is smooth it must be tangent to the
distribution ${\cal D}$ normal to ${\rm Sing}(Y_b)$, a
contradiction. This finishes the proof of (3).
 $\Box$

\medskip
{\bf Proposition 4.8} {\it Suppose $Y_b$ is not irreducible. Then
each rational curve $C'$ on $Y_b$ is smooth. Moreover, for each
point $x' \in C'$, the irreducible component of the germ of $Y_b$ at
$x'$ containing the germ of $C'$ at $x'$ is smooth. }

\medskip
{\it Proof}. A rational curve $C'$ on $Y_b$ must be the image of the
Albanese fiber $C$ of some component $\hat{X}$ by Proposition 3.8.
Thus, the statement follows if $\nu_C : C \longrightarrow C'$ is an
isomorphism. Note that $\nu_C$ is of degree $1$, as $C$ is not
contained in $\nu^{-1}({\rm Sing}\, (X))$ by Proposition 3.6.

Suppose that $\nu_C$ is not an isomorphism for some $C$. Since
$\nu_C$ is of degree $1$, the image $C' = \nu(C)$ must be singular
at some point, say at $x' \in C'$. If $X$ is locally irreducible at
$x'$, we get contradiction from Proposition 4.7. Thus $X$ is not locally
irreducible at $x'$ and each component of the germ at $x'$ is smooth
by Proposition 4.1.

Hence the germ of $C'$ at $x$ has at least two irreducible components.
Therefore, $\nu_{C}^{-1}(x')$ contains at least two different
points, say $x_{1}$ and $x_{2}$.

The twisted vector field $0 \not= \hat{\lambda}_H\vert_C \in
H^{0}(C, T(C) \otimes \nu^* {\mathcal O}(H))$ vanishes at $x_{1}$
and $x_{2}$ by Proposition 3.3. On the other hand,  since $Y_b$ is
not irreducible, we have ${\rm deg}\, T(C) \otimes \nu^*{\mathcal
O}(H) \le 1$ by Proposition 4.6, a contradiction. $\Box$

\medskip
For the statement of the next Proposition, it is convenient to
introduce the notion of the irreducible germ of $Y_b$ along a
rational curve. Suppose $Y_b$ is not irreducible. Let $C'$ be a
rational curve on $Y_b$, which is necessarily smooth by Proposition
4.8. Since $C'$ is smooth, at each point $x' \in C'$, there is a
unique irreducible component, say ${\mathcal V}_{x'}$, of the germ
of $Y_b$ at $x'$ containing the germ of $C'$ at $x'$. By Proposition
4.8, ${\mathcal V}_{x'}$ is smooth. Other irreducible components of
the germ of $Y_b$ at $x'$ intersect properly with $C'$. Since $C'$
is compact and irreducible, we can choose $\{ {\mathcal V}_{x'}, \;
x' \in C'\}$ so that $\cup_{x' \in C'} {\mathcal V}_{x'}$ forms a
germ of a complex submanifold in $M$ containing  $C'$. We call
$\cup_{x' \in C'} {\mathcal V}_{x'}$ the {\it irreducible germ of
$Y_b$ along $C'$}.

\medskip
{\bf Proposition 4.9 } {\it Suppose $Y_b$ is not irreducible. Let
$C'$ be a (necessarily smooth) rational curve on $Y_b$. Let
${\mathcal X}_0$ be the irreducible germ of $Y_b$ along $C'$. Let
${\mathcal X}_1, \ldots , {\mathcal X}_{\ell}$ be the other
components of the germ of $Y_b$ meeting $C'$ at some points, and let
${\mathcal H}_i$ be the unique irreducible component of the germ $Y
= f^{-1}(D)_{\rm red}$ at $C' \cap {\mathcal X}_i$ such that
${\mathcal X}_i \subset {\mathcal H}_i$. Denote by $a_i$ the
multiplicity of ${\mathcal H}_i$ in $f^{-1}(D)$. Then,
$$2a_0 = a_1 (C'
\cdot {\mathcal H}_1) + \cdots +
a_{\ell} (C' \cdot {\mathcal H}_{\ell})\,\,.$$}

\medskip
Note that $a_i$ coincides with the multiplicity of the global
irreducible component $H$ of $f^{-1}(D)$ such that ${\mathcal H}_i
\subset H$. The intersection number $(C' \cdot {\mathcal H}_{i})$ is
the one counted with multiplicities at the intersection point. The
intersection number $(C' \cdot {\mathcal H}_{0})$ also makes a
sense, as ${\mathcal O}({\mathcal H}_{0})$ can be naturally regarded
as a line bundle on the complete curve $C'$.

\medskip
{\it Proof}. Let $h$ be the defining equation of the critical
divisor $D$ at $b$. Let $x' \in C'$ and let $h_i$ be the local
equation of ${\mathcal H}_i$ at $x'$. Then $f^* h= h_0^{a_0} \cdot
h_1^{a_1} \cdots h_{\ell}^{a_{\ell}}$ in ${\mathcal O}_{M, x'}$.
Note that ${\mathcal H}_0$ is a globally defined divisor around
$C'$. Thus, we can consider the ${\mathcal O}(-(a_0 -1){\mathcal
H}_0)$-valued vector field $\gamma$ in a neighborhood of $C'$ in
$M$, which is locally defined by
$$\gamma = \frac{1}{h_0^{a_0-1}}\iota_{\omega}(df^* h) = a_0 h_1^{a_{1}} \cdots
h_{\ell}^{a_{\ell}} \iota_{\omega}(d h_0) + \gamma_1  \,\, ,$$ where
$\gamma_1$ is a vector field vanishing on ${\mathcal H}_0$. Note
that the  adjunction formula gives
$${\mathcal H}_{0} \cdot C' = K_{{\mathcal H}_0} \cdot C' =
{\rm deg}\, K_{C'} = -2\,\, $$ because $K_M$ and the normal bundle
of $C'$ in ${\mathcal H}_0$ are both trivial. By the same argument
as in the proof of Proposition 2.1, $\gamma$ is tangent to $C'$
and it gives rise to a non-zero global section $\gamma \vert_{C'}$
of the line bundle $T(C') \otimes {\mathcal O}(-(a_0 -1){\mathcal
H}_0)$ on $C'$. Since $C' = {\mathbf P}_1$, this line bundle is of
degree
$$-(a_0-1) (C' \cdot {\mathcal H}_0) +2 = 2(a_0 -1) + 2 = 2a_0\,\, .$$
On the other hand, from the local expression of $\gamma$ above and
the fact that ${\mathcal H}_0$ is smooth along $C'$, we see that
$\gamma \vert_{C'}$ has exactly
$$a_1 (C' \cdot {\mathcal H}_1) + \ldots + a_{l} (C' \cdot {\mathcal H}_{\ell})$$
zeros counted with multiplicities.
It follows that
$\sum_{i =1}^{\ell} a_i (C' \cdot {\mathcal H}_i) = 2 a_0$. $\Box$

\medskip
Now we are ready to prove Theorem 1.4 when $Y_b$ is singular. By construction,
the ${\mathbf C}^{n-1}$-action in Proposition 2.2 descends to the action
on the set of characteristic $1$-cycles on $Y_b$. In particular, they are
all isomorphic. There
are the following three  possibilities.

\medskip
(Case i) Each singular germ of $Y_b$ falls into the (Case 3) in
Proposition 4.1 and all rational curves on $Y_b$ are smooth.

(Case ii) Some singular germ of $Y_b$ falls into (Case 1) or (Case
2) in Proposition 4.1 and all rational curves on $Y_b$ is smooth.

(Case iii) $Y_b$ is irreducible and all rational curves on $Y_b$ are
singular.

\medskip
In fact, if $Y_b$ is not irreducible, its germ at any singular
point is reducible, thus falls into (Case i) or (Case ii) by
Proposition 4.8. If $Y_b$ is irreducible and a rational curve on
$Y_b$ is smooth, it belongs to (Case i) or (Case ii) by
Proposition 4.7.

\medskip
{\bf Proposition 4.10} {\it Suppose $Y_b$ is in (Case i). Then the
characteristic $1$-cycle $\sum_{s} \frac{r_s}{r} \Theta_s$ divided
by $r = {\rm GCD}\{r_s\}$ is either one of singular fibers of a
relatively minimal elliptic fibration listed by Kodaira [Kd,
Theorem 6.2], $1$-cycle of Type $A_{\infty}$ or $1$-cycle of Type
$D_{\infty}$, defined in Theorem 1.4.}

\medskip
{\it Proof}. Let $X_1, \ldots, X_N$ be the (global) irreducible
components of $Y_b$ and let $H_i$ be the (global) irreducible
component of $Y = f^{-1}(D)_{\rm red}$ such that $X_i \subset
H_i$. We denote by $C_{i}'$ any smooth rational curve in $X_i$.
Each $\Theta_s$ is ${\mathbf P}_1$ by our assumption, and it
belongs to one of $X_1, \ldots, X_N$. Let us denote by $\{
\Theta^{j}_i, 1 \leq j \leq \tau(i) \}$, the collection of the
components contained in $X_i$ the cardinality  $\tau(i)$ of which
may be infinite. Let us define the intersection number $\Theta^j_i
\cdot \Theta^p_m$ as the cardinality of the set-theoretic
intersection $\Theta^j_i \cap \Theta^p_m$, except when $i=m$ and
$j=p$ in which case we {\it decree} it to be $-2$. Set $a_{i m} :=
C_{i}' \cdot H_{m}$. This number does not depend on the choice of
$C_{i}'$ in $X_i$ by Proposition 3.8. Note that for each $1 \leq i
\not= m \leq N$ and for each $j$, there are only finitely many $p$
such that $\Theta^j_i \cap \Theta_{m}^{p} \not= \emptyset$.
Indeed, in (Case i), such curves $\Theta^p_m$ bijectively
corresponds to the set $\Theta^j_i \cap {\rm Sing}\, (Y_b)$, which
consists of at most finitely many points. Thus, for each $1 \leq i
\neq m \leq N$ and $ 1 \leq i \leq \tau(i),$ we have a
well-defined equality
$$ a_{i m} = \Theta^j_i \cdot H_m = \sum_{p} \Theta^j_i \cdot \Theta^p_m
\,\, .$$
This equation is also true when $i = m$, as
$$a_{ii} = \Theta^j_i \cdot H_i = \Theta^j_i \cdot {\mathcal H}_{i, 0} +
\Theta^j_i \cdot {\mathcal H}_{i, 1} +  \cdots + \Theta^j_i \cdot
{\mathcal H}_{i, q}\,\, {\rm and} \,\, \Theta^j_i \cdot {\mathcal
H}_{i, 0} = -2 = \Theta^j_i \cdot \Theta^j_i\,\, .$$ Here ${\mathcal
H}_{i, 0}$ is a germ of $H_i$ along $\Theta^j_i$ and ${\mathcal
H}_{ij'}$ ($j' \ge 1$) are other germs of $H_i$ which meets
$\Theta^j_i$ at some points.

Then, by $\Theta^j_i \cdot f^{-1}(D) = 0$, we have
$$\sum_{1 \leq m \leq N, 1 \leq p \leq \tau(m)} r_m \Theta^j_i
\cdot \Theta^p_m = \sum_{1\leq m \leq N} r_m a_{i m} = 0\,\, .$$
This means the intersection numbers $\Theta_t \cdot \Theta_s$ and
the coefficient $r_s$ of the cycle $\sum_s r_s \Theta_s$ satisfy $$
2 r_s = \sum_{t \neq s} r_t \Theta_t \cdot \Theta_s\,\, .$$ This
formula, which we will call the key formula later, is exactly
the formula in [Kd, p. 567, equation (6.5)]. The only difference
is that the cardinality of indices here is possibly infinite. Now
all the arguments in [Kd, pp. 567--571] for the classification of
$\sum_s r_s \Theta_s$, starting from a component with minimal multiplicity,
say $r_0$, goes through modulo the following three points:

(a) the great common divisor $r$ of $r_s$ may not be $1$ (even if it is not of
Type $I_n$);

(b) the process of [Kd, pp 568, lines 23--25, the case ($\beta$)] is not terminate if $\sum_s r_s \Theta_s$ is an infinite $1$-cycle, so that we could {\it not} conclude that $\sum_s r_s \Theta_s$ contains a (finite) cyclic chain even if $\Theta_0$ meets at least two components;

(c) the process of [Kd, pp 569, the case ($\beta_2$)] is not terminate if $\sum_s r_s \Theta_s$ is an infinite $1$-cycle.

For (a), we may just divide the cycle by $r$. From the key formula above,
the cycle is of Type $A_{\infty}$ when the process (b) is not terminate.
Again, from the key formula above, the cycle is of Type $D_{\infty}$
when the process (c) is not terminate. This completes the proof. $\Box$

\medskip
{\bf Proposition 4.11} {\it When $Y_b$ is in (Case ii),  one of the
following holds.

(1) If a local germ of $Y_b$ contains two components, there are at
most two components in $Y_b$ and the singular locus of $Y_b$ is
irreducible. The multiplicities of $f^{-1}(D)$ along all
components of $Y_b$ are equal. Moreover through each point of the
singular locus there are exactly two smooth rational curves,
meeting each other with contact order $2$. In this case, the characteristic $1$-cycles are isomorphic to the Kodaira fiber of Type $III$ in [Kd]
up to the common multiplicities. (Both an example of an
irreducible $Y_b$ and an example of $Y_b$ with two components
exist; see e.g. [M2].)

(2) If a local germ of $Y_b$ contains three components, there are
at most three components in $Y_b$ and the singular locus of $Y_b$
is irreducible. The multiplicities of $f^{-1}(D)$ along all
components of $Y_b$ are equal. Moreover through each point of the
singular locus there are exactly three smooth rational curves,
meeting each other transversally. In this case, the characteristic
$1$-cycles are isomorphic to the Kodaira fiber of Type $IV$ in [Kd] up to the
common multiplicities. (There are examples of this type where the
number of components of $Y_b$ can be one, two or three; see e.g.
[M2].) }

\medskip
{\it Proof}. Let $\sum_s r_s \Theta_s$ be a characteristic
$1$-cycle. By our assumption, $\Theta_s = {\mathbf P}_1$.
Define the intersection number $\Theta^i_j \cdot \Theta^k_m$ by
the cardinality of the set-theoretic intersection, weighted by the
order of contact. The order of contact can be 1 or 2 by
Proposition 4.1. We also put $\Theta_s \cdot \Theta_s = -2$ as in
the proof of Proposition 4.10. Then the same argument as in the
proof of Proposition 4.10, following [Kd, proof of Theorem 6.2], goes
through. By the assumption in (Case ii) and by
the shape of singular fibers in Kodaira's list [Kd], our characteristic
$1$-cycle is a finite $1$-cycle and, as $1$-cycles, it would be of
Type $III$ if $Y_b$ had a singular germ in (Case 2) of Proposition 4.1
and of Type $IV$ if
$Y_b$ had a singular germ in (Case 1) of Proposition 4.1.

Since Kodaira fibers of type $III$ and $IV$ are not semi-normal
(see [GPR, Chapter I, Section 15] for a definition), it is not
immediate that our cycle $\sum_s r_s \Theta_s$ is actually {\it
biholomorphic} to the Kodaira fiber. In other words, the equality
of the coefficients and the intersection numbers only guarantees
that the semi-normalization of the fiber of $\tilde{\alpha}$ is
isomorphic to the semi-normalization of the corresponding Kodaira
fiber. The biholomorphicity can be seen as follows.

In the case of Type $III$, since the germ of  two smooth analytic
curves having contact of order 2 is always  biholomorphic to the
germ of two such curves on a smooth surface, it is immediate that
the Albanese fiber is biholomorphic to the Kodaira fiber.

In the case of Type $IV$, it is easy to see that the germ of three
smooth curves intersecting at one point, with transversal pairwise
intersection, is isomorphic to the germ of the Kodaira fiber if and
only if the three tangent vectors at the intersecting point are
linearly dependent. But this is true for the Albanese fiber by
Proposition 4.4 (3). Thus the Albanese fiber is biholomorphic to the
Kodaira fiber.

All the other statements in Proposition 4.11 follow immediately from
the properties of the corresponding Kodaira fibers.
 $\Box$

\medskip
{\bf Proposition 4.12} {\it When $Y_b$ is in (Case iii), one of the
following holds.

(1) Each rational curve on $Y_b$ is isomorphic to the cuspidal
cubic curve in the plane. The singular locus of $Y_b$ is
irreducible and consists of the cusps of the rational curves. In
this case, the characteristic $1$-cycles are irreducible and of
Type $II$ in [Kd] up to multiplicity.

(2) Each rational curve on $Y_b$ is an immersed ${\bf P}_1$ with one
node. In this case, $Y_b$ is semi-normal and the characteristic $1$-cycles are
irreducible and of Type $I_{1}$ in [Kd], up to multiplicity. }

\medskip
{\it Proof}. If an irreducible component of the germ of $Y_b$ is
singular, then we are in the situation of Proposition 4.7 and (1)
follows. If each component of the germ of $Y_b$ is smooth, then a
rational curve on  $Y_b$ is an immersed ${\bf P}_1$ with
self-intersection. Since $Y_b$ is irreducible, $\lambda_H$ is a
vector field along $C$ vanishing at the points over the
self-intersection points as in the proof of Proposition 4.7. The
germ of $Y_b$ at the singular point can have transversal
intersection or contact of order 2. If it has contact of order 2,
we see that the vector field $\hat{\lambda}_H \vert_C$ has double
zeroes at the two points of $C$ over the self-intersection, a
contradiction. Thus there is only one nodal point and we are in
(2). $\Box$

\medskip
As Matsushita pointed out to us, there actually occurs the  case
where the characteristic $1$-cycle is an infinite cycle of Type
$A_{\infty}$. The following proposition and its proof are due to
Matsushita:

\medskip
{\bf Proposition 4.13} {\it There is a $4$-dimensional Lagrangian fibration
$f : M \longrightarrow \Delta^2$ such that $f$ is projective and such
that the characteristic $1$-cycles of a general singular fiber of $f$
are of Type $A_{\infty}$. }

\medskip
{\it Proof} Let $R_{k} = {\rm Specan}\, {\mathbf C}[u^{k+1}v^{-1},
u^{-k}v]$ ($k \in {\mathbf Z}$).  There is a natural morphism
$g_{k} : R_{k} \longrightarrow {\rm Specan}\, {\mathbf C}[u]$. Let
$E$ be an elliptic curve. Using the morphisms $g_{k}$, which are
compatible with the natural gluing of the spaces $R_{k}$, we
obtain a morphism
$$(\cup_{k \in {\mathbf Z}} R_{k}) \times E \times {\rm Specan}\, {\mathbf C}[y] \longrightarrow {\rm Specan}\, {\mathbf C}[u, y]\,\, .$$
Restricting this morphism over a sufficiently small $2$-dimensional disk $\Delta^2$ (centered at $(u, y) = (0,0)$), we obtain a fibration
$$\tilde{f} : \tilde{M} \longrightarrow \Delta^{2}_{(u, y)}\,\, .$$
The fiber over $u = 0$ is an infinite chain of ${\mathbf P}_1 \times E$,
while the fiber over $u \not= 0$ is ${\mathbf C}^{*} \times E$. Let $\alpha$ be a {\it non-torsion} point of $E$. Then
${\mathbf Z}$ acts on $\tilde{f} : \tilde{M} \longrightarrow \Delta^2$ if we define the action of $m \in {\mathbf Z}$ by
$$(u^{k+1}v^{-1}, u^{-k}v, x, y) \mapsto (u^{k+1+m}v^{-1}, u^{-k-m}v, x + m\alpha, y)\,\, .$$
By [Na, Theorem 2.6], this action is properly discontinuous (and is free). Moreover, by [Na, Section 5], the induced morphism
$$f : M = \tilde{M}/{\mathbf Z} \longrightarrow \Delta^{2}$$
is projective. As the symplectic $2$-form $du \wedge dv/v + dx
\wedge dy$ on $\tilde{M}$ is ${\mathbf Z}$-invariant, it descends
to a symplectic $2$-from on $M$. With respect to this form, $f$ is
a Lagrangian fibration. By construction, the normalization of each
singular fiber is ${\mathbf P}_1 \times E$. Characteristic curves
are exactly the image of the fibers ${\mathbf P}_1$ of ${\mathbf
P}_1 \times E$. Since $\alpha$ is of infinite order, the
characteristic $1$-cycles of $f$ are then of Type $A_{\infty}$.
$\Box$

\bigskip
{\bf References}

\medskip

[Ar] Arnold, V. I. : Mathematical methods of classical mechanics.
Graduate Texts in Mathematics {\bf 60}, Springer-Verlag, New York, 1989.

[DM] Donagi, R. and Markman, E.: Spectral covers, algebraically
completely integrable, Hamiltonian systems, and moduli of bundles.
Integrable systems and quantum groups, 1--119, Lecture Notes in
Math., {\bf 1620}, Springer, Berlin, 1996.

[Fu] Fujiki, A.: On automorphism groups of compact K\"ahler
manifolds. Invent. math. {\bf 44} (1978) 225--258.

[GPR] Grauert, H., Peternell, T. and Remmert, R. : Several Complex
Variables VII, Encyclopedia of Mathematical Sciences {\bf 74},
Springer-Verlag 1994.

[Ho] Holmann, H.: Seifertsche Faserraume.  Math. Ann. {\bf 157} (1964)
138--166.

[HR] Hwang, J.-M. and Ramanan, S.: Hecke curves and Hitchin
discriminant. Ann. scient. Ec. Norm. Sup. {\bf 37} (2004) 801--817.

[Kd] Kodaira, K.: On compact analytic surfaces. II,
Ann. of Math. (2) 77 (1963) 563--626.

[M1] Matsushita, D.: Higher direct images of dualizing sheaves of
Lagrangian fibrations,  AG/0010283v1, Amer. J. Math. {\bf 127} (2005)
243--259.

[M2] Matsushita, D.: On singular fibers of Lagrangian fibrations
over holomorphic symplectic manifolds. Math. Ann. {\bf 321} (2001)
755--773.

[M3] Matsushita, D.: A canonical bundle formula for projective
Lagrangian fibrations. preprint, 2007.

 [Na] Nakamura, I.: Relative compactification
of the N\'eron model and its application; In Complex analysis and
algebraic geometry, Iwanami Shoten (1977) 207--225.

[Og] Oguiso, K.: Local families of K3 surfaces and applications, J.
Alg. Geom. {\bf 12} (2003) 405--433.

[Ra] Ramanujam, C. P.: On a certain purity theorem, J. Indian Math. Soc.
{\bf 34} (1970) 1--9.

[Sa] Sawon, J.: Foliations on hypersurfaces in holomorphic
symplectic manifolds, preprint, 2007.

[Se] Serre, J.-P. : Algebraic groups and class fields.
Graduate Texts in Mathematics {\bf 117}, Springer-Verlag, New York, 1988.

[Si] Siu, Y.-T.:  Nondeformability of the complex projective space.
J. Reine Angew. Math. {\bf 399} (1989), 208--219, and Errata. J.
Reine Angew. Math. {\bf 431} (1992), 65--74.

\pagebreak

Jun-Muk Hwang

Korea Institute for Advanced Study

207-43 Cheongnyangni-dong

Seoul 130-722, Korea

 jmhwang@kias.re.kr

\bigskip

Keiji Oguiso

Keio University

4-1-1 Hiyoshi Kouhoku-ku 223-8521,

Yokohama, Japan $\;\;$ {\it and}

Korea Institute for Advanced Study

207-43 Cheongnyangni-dong

Seoul 130-722, Korea

oguiso@hc.cc.keio.ac.jp

\end{document}